\documentclass[12pt]{article}
\usepackage{amssymb}
\usepackage{amsmath}

\usepackage{color}
\usepackage{geometry}

\newif\ifbb
\bbfalse

\begin{document}


\newtheorem{lem}{Lemma}
\newtheorem{thm}{Theorem}
\newtheorem{cor}{Corollary}
\newtheorem{prp}{Proposition}

\newtheorem{rem}{Remark}

\newcommand\tz{\kern.3mm{\rm{;}}\ }
\newcommand\z{{\kern.3mm\rm{,}}\ }
\newcommand\dv{\kern.3mm{\rm{:}}\ }
\renewcommand\t{\kern.3mm{\rm{.}}\ }
\newcommand\vop{\kern.3mm{\rm{?}}\ }
\renewcommand\[{$[$}
\renewcommand\]{$]$}
\newcommand\ts{{\text{\footnotesize(}}}
\newcommand\tss{{\text{\footnotesize)}}}

\newcommand\wtb{{}\,\widetilde{\!B}\kern.1mm}


\title{{\normalsize\tt\hfill\jobname.tex}\\
On partial derivatives \break  of multivariate Bernstein polynomials
\footnote
{For the first author the work was prepared within the framework of a subsidy granted to the HSE by the Government of the Russian Federation for the implementation of the Global Competitiveness Program.}
}
\author{A.~Yu.~Veretennikov\footnote{University of Leeds, UK; National Research University 
``Higher School of Economics''
and Institute for Information Transmission Problems, Moscow, Russia;
e-mail: a.veretennikov @ leeds.ac.uk}, 
E.~V.~Veretennikova\footnote{Moscow State Pedagogical University, Russia.}}

\maketitle

\begin{abstract}
\noindent It is shown that Bernstein polynomials for a multivariate function converge to this function along with partial derivatives provided that the latter derivatives exist and are continuous. This result may be useful in some issues of stochastic calculus. 
\end{abstract}

\section{Introduction}

Widely known is the proof of the polynomial Weierstrass theorem based on Bernstein
polynomials and on the Law of Large Numbers for Bernoulli trials proposed in~\cite{bern}.

This method is applicable for the approximation  of multivariate functions, too. In the literature it was noted that for a univariate function, Bernstein polynomials approximate this function along with its derivatives assuming that they exist, see for example \cite{Hlo, deriv}. It turned out that this observation about an approximation along with derivatives holds true for multivariate   functions as well; however, the authors were unable to find an exact reference and this was the reason for writing the present paper. In  particular, this property is very useful in one of the proofs of multi-dimensional Ito's formula, which formula is the main tool in stochastic analysis  \big(cf. \cite[Theorem 10.4]{Kry}\big).

The latter proof starts
with a verification of Ito's formula for the products; then it is extended by induction from linear functions to arbitrary polynomials and finally to all functions
with two continuous and bounded derivatives; this approach essentially follows the lines of one of the proofs of the Weierstrass theorem \big(cf.~\cite[ch.\,XVI, \S\,4]{Zor2}\big).
It is essential that not only the function
itself, but also its derivatives up to the second order were approximated on any
compact. The latter requirement makes this version of the Weierstrass theorem a bit non-standard.
At least, the majority
of textbooks \big(cf.~\cite{{Zor2},{ash}}\big) where approximations by polynomials are discussed,  
usually focus on functions themselves and not on their derivatives. 

In the univariate case this nuance is not important because we can approximate the $n$-th order derivative and then integrate it $n$ times. But in the multivariate case 
this trick does not help because integrals may depend on the paths.
In \cite{Kry} the property of approximation of a function along with its derivatives is simply claimed as widely known. The authors found it interesting to inspect whether celebrated Bernstein polynomials, which gave rise to a notable branch of approximation theory \big(cf.~\cite{review}\big) admit this property. 
The paper \cite{VV2012} presents briefly the main result and the idea of approximation of partial derivatives in the case of~$d=2$.
In this paper a full proof of this fact in the general case of $d\ge 2$ is provided. The 
case $d=1$ is briefly presented in the next paragraph for the sake of completeness. 

This paper consists of five sections. The first is Introduction. In the second some classical theorems and well-known  generalizations in one-dimensional case are recalled. The third contains the main results about multivariate Bernstein
polynomial derivatives convergence. The section 4 and 5 present the proofs of an auxiliary lemma  
and of the main results, correspondingly.

\section{The case $d=1$}                                    
For any function $f$ on $[0,1]$, approximating Bernstein polynomials are given by the formula
\vskip-0.00007pt\noindent
$$
B_n(f;x):=\sum_{j=0}^nf
\bigg(\frac{j}{n}
\bigg)C^j_nx^j(1-x)^{n-j},\ \ 0\le x\le 1.
$$

\begin{thm}[Bernstein]\label{thm1-en}
If $f\in C([0,1])$, then $B_n(f;x) \to f^{}(x), \; n\to\infty$ and this convergence is uniform $[0,1]$.
\end{thm}
See \cite{bern}, and some generalisations in \cite{tel}. The following similar result holds true for the derivatives.
\begin{thm}\label{thm2-en}
If $f\in C^k([0,1])$, then $B^{(k)}_n(f;x) \to f^{(k)}(x), \; n\to\infty$ and convergence is uniform on $[0,1]$.
\end{thm}
\goodbreak
The proof of the Theorem \ref{thm2-en} see, e.g., in \cite{L, L1937}; the result was established by I.~N.~Chlodovsky \cite{Hlo}. \big(The latter reference \cite{Hlo} gives only the title of his talk at the All-Unions Mathematical congress; the texts of most talks have not been published. According to \cite[Chapter 1]{Ti-Sh}, since then the result of the  Theorem \ref{thm2-en} became known in the literature. However, the issue of correct references is a bit unclear.\big)
There are various bounds of convergence rate under additional assumptions about smoothness or without them \big(cf. \cite{rate, deriv, pop35} et al.\big), but they are not the goal of this paper.
The following identity for derivative will be useful in the sequel \big(cf., eg., \cite{L, L1937}\big):
\begin{equation}\label{eq1}
B'_n(f;x)=
\sum_{j=0}^{n-1}n\Delta_{1/n}f
\bigg(\frac{j}n
\bigg)C^j_{n-1}x^j(1-x)^{n-1-j},
\end{equation}
where
\begin{equation}\label{Deltaf} 
\Delta_zf\bigg(\frac{j}n
         \bigg):=f
         \bigg(\frac{j}n+z
         \bigg)-f
         \bigg(\frac{j}n
         \bigg),
\quad z\in\mathbb R.
\end{equation}
Later for a multi-dimensional case we will use a more detailed notation $\Delta_{z, x^i}$, which emphasises that the increment corresponds to the  variable $x^i$; for a fixed $z=1/n$ we will use a short notation $\Delta_{(x^i)}$. By induction the representation of the derivative of order $k$ follows:
\begin{equation}\label{eq3}  
B^{(k)}_n(f;x)=
\frac{n!}
     {(n-k)!}\,
\sum_{j=0}^{n-k}\Delta_{1/n}^kf
\bigg(\frac{j}n
\bigg)C^j_{n-k}x^j(1-x)^{n-k-j},
\end{equation}
where $\Delta_{z}^k$ is defined also by induction as the operator $\Delta_{z}$ applied $k$ times. 
For example for $k=2$ we have,
\begin{eqnarray*}
\Delta_z^2f(x)=\Delta_z
\big(\Delta_zf(x)
\big)=f(x+2z)-2f(x+z)+f(x).
\end{eqnarray*}
Note that the latter expression is one of the versions of a non-normalized finite difference  Laplacian for the one-dimensional case.
Further, since $n^k\Delta_{1/n}^k f(x)\rightrightarrows f^{(k)}(x)$ 
$\displaystyle \frac{n!}{(n-k)!}/n^k \to 1, \, n\to\infty$, under any fixed $k$, 
then by virtue of (\ref{eq3}) and due to the main calculus in the proof of the Theorem \ref{thm1-en} based on the law of large numbers in the Bernoulli trials scheme with probability of success $x$ (we recall that $x\in [0,1]$), the proof of  convergence for the derivatives follows immediately, of course, under the assumption $f\in C^{k}(\mathbb R)$. It should be noted that all convergences are uniform on $[0,1]$. 
We do not show here the details of this well-known calculation \big(cf. \cite[Theorem 1.8.1]{L} in the case $k=1$\big) as for the main result of this paper in the multi-dimensional case they all will be given in the next sections. Yet we want to point out the development of the key Bernstein's idea,  which suggests to reformulate
(\ref{eq1}) and (\ref{eq3}), respectively, as 
\begin{align*}
B'_n(f;x)&=\mathbb E n\Delta_{1/n}
f\big(n^{-1}\xi_{n-1}(x)
 \big)\\\noalign{\vskip-12pt}
\intertext{and}\noalign{\vskip-12pt}
B^{(k)}_n(f;x)&=\mathbb E\frac{n!}{(n-k)!}\,\Delta_{1/n}^kf
\big(n^{-1}\xi_{n-k}(x)
\big),
\end{align*}
where $\xi_n(x)$ denotes a random variable with Binomial distribution $\mbox{Bin}(n,x)$.
These representations clearly explain why in the case $d=1$ the left hand sides of (\ref{eq1}) and (\ref{eq3}) tend by the law of large numbers to their limits $f'(x)$ and $f^{(k)}(x)$, respectively, under the  assumptions of the Theorem \ref{thm2-en}. 

\section{The case $d\ge 2$~--- Main Result}\label{Sec3}  

Now let us consider $d\ge 2$, $x\in \mathbb R^d$, $x=(x_1, \dots , x_d)$.
There are at least two ways to define Bernstein polynomials in the multi-dimensional case (actually, there are many ways and we will describe them later): either on a simplex $\mathbb S^d:=((x_1,x_2, \dots , x_d):\, 0\le x_1, x_2, \dots , x_d,\;\; \|x\|
\le 1)$ by the formula
\begin{equation}\label{b} 
B_n(f;x)
:=\hskip-10pt
\sum_{\substack{0\le j_1+j_2+\dots+j_d\le n,\\j_1,j_2,\dots,j_d\ge 0}}
\!\!f\bigg(\frac{j_1}n,
       \frac{j_2}n,\dots,
       \frac{j_d}n
 \bigg)  C^j_nx_1^{j_1}x_2^{j_2}\cdots x_d^{j_d}
\big(1-\|x\|
\big)^{n-|j|},
\end{equation}
where $j=(j_1,j_2, \dots , j_d)$, and where a vector $x$ norm, and the ``modulus'' of a multi-index $j$, and the polynomial coefficient (``$n$ choose $j$'') $C^{j}_{n}$ are defined as
$\|x\| = |x_1|+|x_2| + \dots + |x_d|$, 
$
|j| = j_1+j_2 + \dots + j_d,
$
and
$$
C^{j}_{n} \equiv C^{j_1, \dots, j_d}_{n} = 
\frac{n!}{j_1! j_2! \dots j_d! \big(n-|j|\big)!};
$$
or on a ``$d$-dimensional square'' (cube, etc.) by a formula similar to but different from (\ref{b}),
\begin{align}\label{tildeb}       
\wtb_n(f;x)
&:=\sum_{j_1,j_2,\dots,j_d=0}^nf
\bigg(\frac{j_1}n,
      \frac{j_2}n,\dots,
      \frac{j_d}n
\bigg)C^{j_1}_nC^{j_2}_n\cdots C^{j_d}_n\nonumber\\
\noalign{\vskip7pt}
&\qquad\times x_1^{j_1}(1-x_1)^{n-j_1}
              x_2^{j_2}(1-x_2)^{n-j_2}\cdots
              x_d^{j_d}(1-x_d)^{n-j_d}.
\end{align}
\goodbreak
Note that the degrees of the polynomials $B_n$ and $\tilde B_n$ are different. The formulae (\ref{b}) and (\ref{tildeb}) allow probabilistic representation

\begin{equation}\label{bb}  
B_n(f;x):=\mathbb Ef\big(n^{-1}\eta_n(x)
                    \big)
\end{equation}
and 
\begin{equation}\label{bbb}   
\wtb_n(f;x):=\mathbb Ef\big(n^{-1}\xi_n(x)
                           \big),
\end{equation}
where the distribution of the random vector $\eta_{n}(x)$ of dimension $d\, $ is the projection onto the first $d$ coordinates of the {\em polynomial (multinomial) distribution} of dimension $d+1$ with parameters $n$ (the number of trials) and the vector of ``success probabilities'' $(x_1, \dots, x_{d}, x_{d+1})$, which satisfy conditions $x_i\ge 0$ ($i = 1, \dots, d+1$) and  $\sum_{i = 1}^{d+1} x_i =~1$, while the random vector $\xi_n(x)= (\xi_n^1, \dots, \xi_n^d)$ consists of independent components distributed binomially $\mbox{Bin}(n,x_i)$ each.

Due to the law of large numbers for Bernoulli trials for each coordinate we have, 
\begin{equation}\label{llnxi}        
\frac1n\xi_n(x)\to x,
\ \ n\to\infty,
\end{equation}
both in probability and almost surely. For the sequence ($\eta_n$) the law of large numbers also holds true:
\begin{equation}\label{llneta}   
\frac1n\eta_n(x)\to x,
\ \ n\to\infty,
\end{equation}
in probability and almost surely. The easiest way to show this is to use the law of large numbers for each coordinate, which follows directly from the one-dimensional version of this theorem for the  binomial distribution. 

It apparently does not help too much to analyse partial derivatives where certain series still need to be treated. However, after the differentiation of these series we will get expressions, which admit  representations via expectations of finite differences for the function $f$ with random vector arguments distributed either according to a polynomial law  -- or, more precisely, its projection -- or as a direct product of binomial distributions, which will eventually lead to the desired result.

\begin{rem}     
Other variants are possible, such as a  combination of a ``cube'' and of a simplex for different variables. For example, in case of $d=d_1 + d_2$ and $d_1\ge 2$ one more version of multivariate  Bernstein polynomials has a form, 
\goodbreak
\noindent
\begin{align*}\label{bhat}
\widehat{B}_n(f;x)
&:=\sum_{\substack{0\le j_1+j_2+\dots+ j_{d_1}\le n,\\
                   j_1,j_2,\dots,j_{d_1}\ge 0}}
f\bigg(\frac{j_1}n,\dots,
       \frac{j_{d_1}}n,
       \frac{j_{d_1+1}}n,\dots,
       \frac{j_{d_1+d_2}}n
 \bigg)\nonumber\\\noalign{\vskip-2pt}
&\quad\times C^{j_1,\dots,j_{d_1}}_n
x_1^{j_1}x_2^{j_2}\cdots x_{d_1}^{j_{d_1}}
\Bigg(1-\sum_{i=1}^{d_1}x_i
\Bigg)^{n-j_1-\dots-j_{d_1}}\nonumber\\
&\quad\times C^{j_{d_1{+}1}}_n\cdots
             C^{j_{d_1{+}d_2}}_nx_{d_1{+}1}^{j_{d_1{+}1}}
             (1{-}x_1)^{n{-}j_{d_1{+}1}}\cdots
             x_{d_1{+}d_2}^{j_{d_1{+}d_2}}
             (1{-}x_{d_1{+}d_2})^{n{-}j_{d_1{+}d_2}},
\end{align*}
where $x=(x_1, \dots, x_{d_1}, x_{d_1+1}, \dots, x_{d_1+d_2})$, all $x_i$ are non-negative, $x_1 + \dots+ x_{d_1} \le 1$, and $0\le x_{d_1+1} \le 1, \dots, 0\le x_{d_1+d_2}\le 1$.

It is quite likely that for such  representations  analogous results about convergence of polynomials and their derivatives may be established.
\end{rem}

So, a multivariate analogue of the Theorem \ref{thm1-en} for the polynomials $B_n$ and $\tilde B_n$ may be formulated as follows.

\begin{thm}\label{thm3-en}
If $f\in C(\mathbb R^d)$, then 
\vskip-2pt\noindent 
\begin{alignat*}2
&\wtb_n(f;x)\to f(x),&&\ \ x\in\mathbb K^d,\\
&B_n(f;x)\to f(x),   &&\ \ x\in\mathbb S^d,
\end{alignat*}
as $n\to\infty$. All convergences are uniform on $\mathbb K^d$ and $\mathbb S^d$, respectively.
\end{thm}

Note that any continuous function on $\mathbb S^d$ or on $\mathbb K^d$ can be extended to a continuous function on $\mathbb R^d$ (of course, not uniquely).  
Convergence of both versions towards a continuous function $f$ on the simplex or on the $d$-dimensional square/cube follows from the (multivariate) law of large numbers: 
in the second case it is applied to a sequence of independent and equally distributed random vectors with independent components, while in the first case -- with dependent components. The proof can be easily found in may papers and textbooks \big(cf. e.g., \cite{deriv}\big) and we omit it.

To state the main result, let us introduce the following notations: :
\begin{enumerate}
\item
$m$~--- a natural number;
\item
$k=(k_1,k_2,\dots,k_d)$~--- a multi-index;
\item
$C^k$ {\em with a multi-index $k$ $k$\/} denotes the class of functions that possess a mixed partial  derivative of order $k= (k_1, k_2, \dots, k_d)$, which is continuous. 
\end{enumerate}
Recall that 
$C^m(\mathbb R^d)$ \big($C^m_b(\mathbb R^d)$\big) with $m = 1, 2, \dots $ is a class of functions on $\mathbb R^d$ with {\em all\/} well-defined mixed derivatives of order $m$, which are continuous (respectively, continuous and bounded).  The notation $f^{(m)}$ stands for the set of all mixed derivatives of the function $f$ of order $m$.

\begin{thm}\label{thm4-en}
$1$. If $m>0$,  $f\in C^m(\mathbb R^d)$, then 
\begin{align*}
&\wtb^{(m)}_n(f;x)\to f^{(m)}(x),\ \ x\in\mathbb K^d,\\
&B^{(m)}_n(f;x)\to f^{(m)}(x),\ \ x\in\mathbb S^d,
\end{align*} 
as $n\to\infty$, and all convergences are uniform on $\mathbb K^d$ and $\mathbb S^d$, respectively. 

$2$. If $k = (k_1, k_2, \dots, k_d)$ is a multi-index and $f\in C^{k}(\mathbb R^d)$, then
\begin{alignat*}2
&\wtb^{(k)}_n(f;x)\to f^{(k)}(x),&&\ \ x\in\mathbb K^d,\\
&B^{(k)}_n(f;x)\to f^{(k)}(x),&&\ \ x\in\mathbb S^d,
\end{alignat*} 
as $n\to\infty$, and all convergences are uniform on $\mathbb K^d$ and $\mathbb S^d$, respectively.
\end{thm}

It is worth mentioning also the papers \cite{ia} and \cite{B2} 
devoted to asymptotic expansions and Taylor's expansions for Bernstein polynomials of two and many variables. While being conceptually close, the results of the present paper do not follow directly from these papers, and we do not aim to get neither asymptotic nor Taylor's expansions here.

\section{Auxiliary result}                      
The following Lemma will be used in the proof of the main Theorem when a uniform convergence of finite differences just under the condition of existence and continuity of the limit expression is needed. While elementary, this Lemma is required for the correctness of references and for a completeness of our presentation.

Let $(z_1, \dots, z_d) \in \mathbb R^d$, and $(x_1, \dots, x_d) \in \mathbb R^d$. Let
\begin{gather*}
\aligned
\Delta_{1/n,x_i}f(x)&:=f(x_1,\dots,x_{i-1},\,x_i+\frac1n,\,x_{i+1},\dots,x_d)\\
&\qquad-f(x_1,\dots,x_{i-1},x_i,x_{i+1},\dots,x_d),
\endaligned\\
\partial^kf(x_1,\dots,x_d):=
\frac{\partial^{|k|}}
     {\partial_{x_1}^{k_1}\dots\partial_{x_d}^{k_d}}f(x_1,\dots,x_d).
\end{gather*}

\begin{lem}\label{lem1}                                    For $f\in C^{k}(\mathbb R^d)$ with $k=(k_1, \dots , k_d)$ and $x=(x_1, \dots , x_d)$ 
the following equality is valid
\begin{gather*}
\Delta^{k_1}_{z_1,x_1}\dots\Delta^{k_d}_{z_d,x_d}f(x)
=\int_{x_1}^{x_1+z_1}d\xi^1_1
  \int_{\xi^1_1}^{\xi^1_1+z_1}d\xi^1_2
\dots\int_{\xi^1_{k_1-1}}^{\xi^1_{k_1-1}+z_1}d\xi^1_{k_1}\\
\dots\int_{x_d}^{x_d+z_d}d\xi^d_1
\int_{\xi^d_1}^{\xi^d_1+z_d}d\xi^d_2\dots
            \int_{\xi^d_{k_d-1}}^{\xi^d_{k_d-1}+z_d}d\xi^d_{k_d}
            \partial^kf(\xi^1_{k_1},\dots,\xi^d_{k_d}).
\end{gather*}
\end{lem}

For example, for $k_1=0$, integration over $\xi^1_1$, $\xi^1_2$, etc. is just not performed, and $f(\xi^1_{k_1},\xi^2_{k_2}, \dots, \xi^d_{k_d})$ is treated as $f(x_{1}, \xi^2_{k_2},\dots, \xi^d_{k_d})$. In particular, if all $k_i=0$ then the equality turns into the identity $f(x)=f(x)$.

The equality for two variables $x=(x_1,x_2)$ and $k_1=k_2=1$
takes the form (for the sake of simplicity we omit the lower indices in $\xi^1_{1},\xi^2_{1}$): 
\begin{equation*}
\Delta_{z_1,x_1}
\Delta_{z_2,x_2}f(x)=
\int_{x_1}^{x_1+z_1}d\xi^1
\int_{x_2}^{x_2+z_2}d\xi^2
\frac{\partial^2}
     {\partial x_1\partial x_2}
f(\xi^1,\xi^2).
\end{equation*}

The proof of the Lemma follows straightforward from the (one-dimensional) first theorem of the calculus (also known as Newton-Leibniz formula) by induction.

\section{Proof of the Theorem~\ref{thm4-en}}        
{\bf 0.} 
It suffices to prove only the second part of the Theorem \ref{thm4-en} because the first part follows immediately due to the identity 
$$
C^m(\mathbb R^d)=\bigcap_{k:\;|k|\le m}C^k(\mathbb R^d).
$$

{\bf 1.} 
Analogues of the formulae~(1) and~(2) in the multivariate case for $\wtb_n$ for the multi-index $k = (k_1, \dots, k_d)$ could be written as
\begin{align}\label{eq33}      
&\frac{\partial^{|k|}}
      {\partial x^k}\wtb_n(f;x)
\equiv\frac{\partial^{|k|}}
           {\partial x_1^{k_1}\dots
            \partial x_d^{k_d}}\wtb_n(f;x)\nonumber\\
\noalign{\vskip7pt}
&\quad=\left(\frac{n\dots(n-k_1+1)\times n\dots(n-k_2+1)
                                  \times\dots\times n\dots(n-k_d+1)}
                  {n^{|k|}}
       \right)\nonumber\\
\noalign{\vskip7pt}
&\qquad
\times\sum_{\substack{0\le j_i\le n-k_i,\\i=1,2,\dots,d}}
 n^{|k|}\Delta^kf\bigg(\frac{j_1}n,
                       \frac{j_2}n,\dots,
                       \frac{j_d}n
                 \bigg)C^{j_1}_{n-k_1}C^{j_2}_{n-k_2}\cdots
                       C^{j_d}_{n-k_d}\nonumber\\
\noalign{\vskip7pt}
&\qquad
\times x_1^{j_1}(1-x_1)^{n-k_1-j_1}x_2^{j_2}(1-x_2)^{n-k_2-j_2}
 \cdots x_d^{j_d}(1-x_d)^{n-k_d-j_d},\qquad
\end{align}
where
$$
\frac{\partial^{|k|}}
     {\partial x^k}:=
\frac{\partial^{|k|}}
     {\partial x_1^{k_1}\dots\partial x_d^{k_d}},
\quad\Delta^k:=\prod_{i=1}^d\Delta_{(x_i)}^{k_i},
$$
and where~$\Delta_{(x_i)}$ means $\Delta_{(x_i)}\equiv\Delta_{1/n,x_i}$.

The usage of the notation $\Delta^k$ with a multi-index $k$ is correct because all operators $\Delta_{(x_i)}$ commute.  Indeed, let us denote (assume $i<j$ for the  sake of definiteness)
\begin{gather*}
x^{i,j}=\big(x^1,\dots,x^{i-1},
             x^i+1/n,x^{i+1},\dots,x^{j-1},
             x^j+1/n,\dots,x^d
        \big),\\
x^i=\big(x^1,\dots,x^{i-1},
         x^i+1/n,
         x^{i+1},\dots,x^d
    \big),\\
x^j=\big(x^1,\dots,x^{j-1},
         x^j+1/n,\dots,x^d
    \big).
\end{gather*}
For $i\not=j$  we get elementary  identities:
\begin{align*}
\Delta_{(x_i)}\big(\Delta_{(x_j)}f(x)
              \big)
&=\Delta_{(x_i)}\big(f(x^{j})-f(x^{})
                \big)\\
&=f(x^{i,j})-f(x^i)-\big(f(x^j)-f(x)
                    \big)\\
&=f(x^{i,j})-f(x^i)-f(x^j)+f(x)
\end{align*}
and similarly, 
\begin{align*}
\Delta_{(x_j)}\big(\Delta_{(x_i)}f(x)
              \big)
&=\Delta_{(x_j)}\big(f(x^i)-f(x)
                \big)\\
&=f(x^{i,j})-f(x^j)-\big(f(x^i)-f(x)
                    \big)\\
&=f(x^{i,j})-f(x^i)-f(x^j)+f(x).
\end{align*}
Analogously to the representation (\ref{bb}), the equality (\ref{eq33}) admits the following probabilistic meaning:
\begin{align}\label{eq333}    
\frac{\partial^{|k|}}
     {\partial x^k}\wtb_n(f;x)
&=\left(\frac{n\dots(n{-}k_1{+}1)\times n\dots(n{-}k_2{+}1)
                                 \times\dots
                                 \times n\dots(n{-}k_d{+}1)}
             {n^{|k|}}
  \right)\nonumber\\
\noalign{\vskip7pt}
&\qquad\times\mathbb En^{|k|}\Delta^kf(n^{-1}\xi_{n-k}(x));
\end{align}
the definition of the random vector $\xi_n(x)=(\xi_n^1,\dots,\xi_n^d)$
see in~\S\,\ref{Sec3}).

{\bf 2.} 
Let us prove the equality~(\ref{eq33}) by induction.
For $k=(0, \cdots , 0)$ (a multi-index of order $d$) the desired equality (\ref{eq3}) is equivalent to the definition of the polynomial $\tilde B^{}_n(f;x)$; this will serve as the {\em basis of induction}. It should be noted that for several variables ($d$ in our case) induction can be carried in turn for the first variable up to $k_1$, then for the second one up to $k_2$, and so on.

In other words, ``double'' induction (it can be also named multivariate) over each of the variables $x_i$, and then over the indices $i=1,\dots,d$ may be applied in our situation. For the induction step, it actually suffices to verify that the formula (\ref{eq3}) remains valid when any component of the multi-index  $k = (k_1,\cdots,k_d)$ increases by one.
For the sake of definiteness let us check the step $k_i \mapsto k_i+1$ for $i=1$.  Recall that 
\begin{align*}
\wtb_n(f;x)
&:=\sum_{j_1,j_2,\dots,j_d=0}^nf
\left(\frac{j_1}n,
      \frac{j_2}n,\dots,
      \frac{j_d}n
\right)C^{j_1}_nC^{j_2}_n\cdots C^{j_d}_n\\
&\qquad\times x_1^{j_1}(1-x_1)^{n-j_1}x_2^{j_2}(1-x_2)^{n-j_2}\cdots
              x_d^{j_d}(1-x_d)^{n-j_d}.
\end{align*}

Because of a certain clumsiness of the formulae and for the sake of clarity we shall start with $k_i=0$ and then proceed to 
the general case. (The same approach will be used for the polynomials $B_n$).
We get, 
\vskip-7pt\noindent
\begin{align*}
\partial_{x_1}\wtb_n(f;x)
&=\partial_{x_1}
  \sum_{j_1,j_2,\dots,j_d=0}^nf
  \left(\frac{j_1}{n},
        \frac{j_2}{n},\dots,
        \frac{j_d}{n}
  \right)C^{j_1}_nC^{j_2}_n\cdots C^{j_d}_n\\\noalign{\vskip2pt}
&\qquad\quad\times x_1^{j_1}(1{-}x_1)^{n{-}j_1}
                   x_2^{j_2}(1{-}x_2)^{n{-}j_2}\cdots
                   x_d^{j_d}(1{-}x_d)^{n{-}j_d}\\\noalign{\vskip2pt}
&=\sum_{\substack{0<j_1\le n,\\0\le j_2,\dots,j_d\le n}}f
  \left(\frac{j_1}n,
        \frac{j_2}n,\dots,
        \frac{j_d}n
  \right)C^{j_1}_nC^{j_2}_n\cdots C^{j_d}_n\\\noalign{\vskip2pt}
&\qquad\quad\times j_1x_1^{j_1{-}1}(1{-}x_1)^{n{-}j_1}
                      x_2^{j_2}(1{-}x_2)^{n{-}j_2}\cdots
                      x_d^{j_d}(1{-}x_d)^{n{-}j_d}\\\noalign{\vskip2pt}
&\qquad-\sum_{\substack{0\le j_1<n,\\0\le j_2,\dots, j_d\le n}}f
\left(\frac{j_1}n, \frac{j_2}n,\dots,
      \frac{j_d}n
\right)C^{j_1}_n
       C^{j_2}_n\cdots
       C^{j_d}_n\\\noalign{\vskip2pt}
&\qquad\quad\times x_1^{j_1}(n{-}j_1)(1{-}x_1)^{n{-}j_1{-}1}
              x_2^{j_2}(1{-}x_2)^{n{-}j_2}\cdots
              x_d^{j_d}(1{-}x_d)^{n{-}j_d}\\\noalign{\vskip2pt}
&=\sum_{\substack{0\le j_1<n,\\0\le j_2,\dots,j_d\le n}}f
  \left(\frac{j_1+1}n,
        \frac{j_2}n,\dots,
        \frac{j_d}n
  \right)C^{j_1+1}_nC^{j_2}_n\cdots C^{j_d}_n\\\noalign{\vskip2pt}
&\qquad\quad\times(j_1+1)x_1^{j_1}(1{-}x_1)^{n{-}(j_1+1)}x_2^{j_2}(1{-}x_2)^{n{-}j_2}\cdots
                    x_d^{j_d}(1{-}x_d)^{n{-}j_d}\\\noalign{\vskip2pt}
&\qquad-\sum_{\substack{0\le j_1<n,\\0\le j_2,\dots,j_d\le n}}f
        \left(\frac{j_1}n,
              \frac{j_2}n,\cdots,
              \frac{j_d}n
        \right)C^{j_1}_n
               C^{j_2}_n\cdots
               C^{j_d}_n\\\noalign{\vskip2pt}
&\qquad\quad\times x_1^{j_1}(n{-}j_1)(1{-}x_1)^{n{-}j_1{-}1}
                   x_2^{j_2}(1{-}x_2)^{n{-}j_2}\cdots
                   x_d^{j_d}(1{-}x_d)^{n{-}j_d}\\\noalign{\vskip2pt}
&=\sum_{\substack{0\le j_1<n,\\0\le j_2,\dots,j_d\le n}}n\Delta_{(x_1)}f
  \left(\frac{j_1}n,
        \frac{j_2}n,\dots,
        \frac{j_d}n
  \right)C^{j_1}_{n-1}
         C^{j_2}_n\cdots
         C^{j_d}_n\\\noalign{\vskip2pt}
&\qquad\quad\times x_1^{j_1}(1{-}x_1)^{n{-}1{-}j_1}
                   x_2^{j_2}(1{-}x_2)^{n{-}j_2}\cdots
                   x_d^{j_d}(1{-}x_d)^{n{-}j_d}, 
\quad
\end{align*}
due to the identities 
\begin{gather*}
(j_1+1)C^{j_1+1}_n=
\frac{(j_1+1)n!}
     {(j_1+1)!\big(n-(j_1+1)
              \big)!}=
\frac{n!}
     {j_1!(n-1-j_1)!}=nC_{n-1}^{j_1},\\
(n-j_1)C^{j_1}_n=
\frac{(n-j_1)n!}
     {j_1!(n-j_1)!}=
\frac{n!}
     {j_1!(n-1-j_1)!}=nC_{n-1}^{j_1}.
\end{gather*}

Now, assuming that the formula (\ref{eq3}) holds true for some $k=(k_1,\dots,k_d)$, let us differentiate it once more in  variable $x_1$ so as to get a similar  formula for the multi-index $(k_1+1,k_2, \dots,k_d)$. For the sake of brevity let us  denote
\begin{align*}
\alpha_{k_1,k_2,\dots,k_d}^{j_1,j_2,\dots,j_d}
&:=C^{j_1}_{n-k_1}\cdots C^{j_d}_{n-k_d}
   x_1^{j_1}(1-x_1)^{n-k_1-j_1}\cdots
   x_d^{j_d}(1-x_d)^{n-k_d-j_d},\\
\noalign{\vskip3pt}
\alpha_{k_2,\dots,k_d}^{j_2,\dots,j_d}
&:=C^{j_2}_{n-k_2}\cdots C^{j_d}_{n-k_d}x_2^{j_2}(1-x_2)^{n-k_2-j_2}
                  \cdots x_d^{j_d} (1-x_d)^{n-k_d-j_d}.
\end{align*}

Note that since we are interested in the  statement of the Theorem for $n\to\infty$,  we can assume that all $k_i \le  n$ \big(although, for $k_i>n$ we have $\partial_{x_i}^{k_i} \tilde B_n(x) = 0$, and the right side of (\ref{eq3}) also equals zero by definition of the number of combinations with negative $n-k_i$\big).

Derivative of a constant being equal to zero, we have, 
\vskip-4pt\noindent
\begin{align*}
&\frac{\partial}
      {\partial x_1}
 \left(\sum_{\substack{0\le j_i\le n{-}k_i,\\i=1,\dots,d}}
       \Delta^kf\bigg(\frac{j_1}n,
                      \frac{j_2}n,\dots,
                      \frac{j_d}n
                \bigg)\alpha^{j_1,\dots,j_d}_{k_1,\dots,k_d}
 \right)\\
&\quad
=\frac{\partial}
       {\partial x_1}
  \left(\sum_{\substack{0\le j_i\le n{-}k_i,\\i=1,\dots,d}}
        \Delta^kf\bigg(\frac{j_1}n,
                       \frac{j_2}n,\dots,
                       \frac{j_d}n
                 \bigg)C^{j_1}_{n{-}k_1}\,x_1^{j_1}(1{-}x_1)^{n{-}k_1{-}j_1}
                 \alpha_{k_2,\dots,k_d}^{j_2,\dots,j_d}
  \right)\\
&\quad
=\frac{\partial}
       {\partial x_1}
  \sum_{j_1=0}^{n{-}k_1}\,
  \sum_{\substack{0\le j_i\le n{-}k_i,\\i=2,\dots,d}}
  \Delta^kf\bigg(\frac{j_1}n,
                 \frac{j_2}n,\dots,
                 \frac{j_d}n
           \bigg)C^{j_1}_{n{-}k_1}\,x_1^{j_1}(1{-}x_1)^{n{-}k_1{-}j_1}
                 \alpha_{k_2,\dots,k_d}^{j_2,\dots,j_d}\\
&\quad
=\sum_{j_1 =1}^{n{-}k_1}\,
 \sum_{\substack{0\le j_i\le n{-}k_i,\\i=2,\dots,d}}
 \Delta^kf\bigg(\frac{j_1}n,
                \frac{j_2}n,\dots,
                \frac{j_d}n
           \bigg)C^{j_1}_{n{-}k_1}\,j_1x_1^{j_1{-}1}(1{-}x_1)^{n{-}k_1{-}j_1}
                 \alpha_{k_2,\dots,k_d}^{j_2,\dots,j_d}\\
&\qquad-\sum_{j_1=0}^{n-k_1-1}
        \sum_{\substack{0\le j_i\le n{-}k_i,\\i=2,\dots,d}}
        \Delta^kf\bigg(\frac{j_1}n,
                       \frac{j_2}n,\dots,
                       \frac{j_d}n
                 \bigg)\\
&\qquad\qquad\times
 C^{j_1}_{n{-}k_1}\,x_1^{j_1}(n{-}k_1{-}j_1)(1{-}x_1)^{n{-}k_1{-}j_1{-}1}
        \alpha_{k_2,\dots,k_d}^{j_2,\dots,j_d}\\
&\quad=\sum_{\substack{0\le j_i\le n{-}k_i,\\i=2,\dots,d}}
       \sum_{j'_1=0}^{n{-}k_1{-}1}\Delta^kf
  \bigg(\frac{j'_1{+}1}n,
        \frac{j_2}n,\dots,
        \frac{j_d}n
  \bigg)\\
&\qquad\qquad\times C^{j'_1{+}1}_{n{-}k_1}\,(j'_1{+}1)
        x_1^{j'_1}(1{-}x_1)^{n{-}k_1{-}j'_1{-}1}
        \alpha_{k_2,\dots,k_d}^{j_2,\dots,j_d}\\
&\qquad-\sum_{\substack{0\le j_i\le n{-}k_i,\\i=2,\dots,d}}
        \sum_{j_1=0}^{n{-}k_1{-}1}\Delta^kf
        \bigg(\frac{j_1}n,
              \frac{j_2}n,\dots,
              \frac{j_d}n
        \bigg)\\
&\qquad\qquad\times
 C^{j_1}_{n{-}k_1}\,x_1^{j_1}(n{-}k_1{-}j_1)(1{-}x_1)^{n{-}k_1{-}j_1{-}1}
              \alpha_{k_2,\dots,k_d}^{j_2,\dots,j_d}\\
&\quad=\sum_{\substack{0\le j_i\le n{-}k_i,\\i=2,\dots,d}}
             \alpha_{k_2,\dots,k_d}^{j_2,\dots,j_d}
  \Bigg(\sum_{j_1=0}^{n{-}k_1{-}1}\Delta^kf
        \bigg(\frac{j_1{+}1}n,
              \frac{j_2}n,\dots,
              \frac{j_d}n
        \bigg)\\
\noalign{\vskip-7pt}
&\qquad\hskip3cm
\times C^{j_1{+}1}_{n{-}k_1}(j_1{+}1)x_1^{j_1}
(1{-}x_1)^{n{-}k_1{-}j_1{-}1}\!{-}\Delta^kf\bigg(\frac{j_1}n,
                \frac{j_2}n,\dots,
                \frac{j_d}n
          \bigg)\\
&\qquad\hskip3cm
\times C^{j_1}_{n{-}k_1}\,
          x_1^{j_1}(n{-}k_1{-}j_1)(1{-}x_1)^{n{-}k_1{-}j_1{-}1}
\Bigg)\\
&\quad=\sum_{\substack{0\le j_i\le n{-}k_i,\\i=2,\dots,d}}
        \alpha_{k_2,\dots,k_d}^{j_2,\dots,j_d}
\left(\vphantom{\sum_{j^{n{-}k_1{-}1}}}
      \sum_{j_1=0}^{n{-}k_1{-}1}(n{-}k_1)\,C^{j_1}_{n{-}(k_1{+}1)}
        x_1^{j_1}(1{-}x_1)^{n{-}(k_1{+}1){-}j_1}
\right.\\
\noalign{\vskip-7pt}
&\hskip1.5cm\qquad\qquad\qquad\times
\left.\vphantom{\sum_{j_{1_1}=0}^{k_1}}
      \Bigg(\Delta^kf\bigg(\frac{j_1{+}1}n,
                           \frac{j_2}n,\dots,
                           \frac{j_d}n
                     \bigg){-}\Delta^kf
                     \bigg(\frac{j_1}n,
                           \frac{j_2}n,\dots,
                           \frac{j_d}n
                     \bigg)\!
      \Bigg)\!
\right)\\
&\quad\equiv\sum_{\substack{0\le j_i\le n{-}k_i,\\i=2,\dots,d}}
            \alpha_{k_2,\dots,k_d}^{j_2,\dots,j_d}
       \sum_{j_1=0}^{n{-}k_1{-}1}(n{-}k_1)\\
\noalign{\vskip-7pt}
&\qquad\qquad\times
 C^{j_1}_{n{-}(k_1{+}1)}x_1^{j_1}(1{-}x_1)^{n{-}(k_1{+}1){-}j_1}
 \Delta_{(x_1)}
 \Delta^kf\bigg(\frac{j_1}n,
                \frac{j_2}n,\dots,
                \frac{j_d}n
          \bigg)\\
&\quad=(n{-}k_1)\sum_{\substack{0\le j_i\le n{-}k_i,\\i=2,\dots,d}}
         \sum_{j_1 =0}^{n{-}(k_1{+}1)}
         \alpha_{k_1{+}1,k_2,\dots,k_d}^{i_1,j_2,\dots,j_d}
         \Delta_{(x_1)}
         \Delta^kf\bigg(\frac{j_1}n,
                        \frac{j_2}n,\dots,
                        \frac{j_d}n
                  \bigg),
\end{align*}
as required. We have used the identities
\begin{align*}
&\aligned
C^{j_1+1}_{n-k_1}(j_1+1)
&=\frac{(n-k_1)!(j_1+1)}
       {(j_1+1)!(n-k_1-j_1-1)!}\\
&=\frac{(n-k_1)(n-k_1-1)!}
       {j_1!(n-k_1-j_1-1)!}=(n-k_1)\,C^{j_1}_{n-(k_1+1)},
\endaligned\\
&C^{j_1}_{n-k_1}\,(n-k_1-j_1)=\frac{(n-k_1-j_1)(n-k_1)!}
                                   {j_1!(n-k_1-j_1)!}
=(n-k_1)\,C^{j_1}_{n-(k_1+1)}.
\end{align*}
Hence, it follows by induction that the formula~(\ref{eq3}) holds true.

{\bf 3.} 
As for $\wtb_n$, it remains to note that for $f\in C^{k}(\mathbb R^d)$ due to the Lemma~\ref{lem1}, appropriately normalized finite differences are uniformly close to the corresponding partial derivatives, that is,
$$
n^{|k|}\Delta^kf(x)
\rightrightarrows
\frac{\partial^{|k|}f}
     {\partial x^k}(x)
\ \ \mbox{���}\ \ x\in\mathbb K^d,\ \ n\to\infty.
$$

So, the statement of the Theorem for the  polynomials $\wtb_n$ follows from (\ref{eq333}) and from the law of large numbers (\ref{llnxi}).

{\bf 4.} 
The analogue of the formula~(\ref{eq3}) for the  polynomials $B_n$ has a form,
\begin{align}\label{eq4}       
\hskip-1mm
B^{(k)}_n(f;x)
&\equiv\frac{\partial^{|k|}}
            {\partial x^k}B_n(f;x)
 \equiv\frac{\partial^{|k|}}
            {\partial x_1^{k_1}\partial x_2^{k_2}\dots
 \partial x_d^{k_d}}B_n(f;x)\nonumber\\\noalign{\vskip2pt}
&=\sum_{0\le j_1+j_2+\cdots\le n-|k|}
  \frac{n!}{\big(n-|k|
            \big)!}\,\Delta_{(x_1)}^{k_1}
                     \Delta_{(x_2)}^{k_2}\dots
                     \Delta_{(x_d)}^{k_d}\nonumber\\\noalign{\vskip4pt}
&\quad\times f\!
\left(\!\frac{j_1}n,
      \frac{j_2}n,\dots,
      \frac{j_d}n\!
\right)C_{n-|k|}^{j_1,j_2,\dots,j_d}x_1^{j_1}x_2^{j_2}\cdots
                                    x_d^{j_d}\big(1{-}\|x\|
                                             \big)^{n{-}|k|{-}|j|}.
\end{align}
For the seequel, let us recall the following notation,
$$
\Delta_{(x_1)}^{k_1}
\Delta_{(x_2)}^{k_2}\dots
\Delta_{(x_d)}^{k_d}=\Delta^k.
$$

The formula~(\ref{eq4}) also has a simple probabilistic meaning. Similarly to the  representation~(\ref{bbb}) we write, 
\begin{align}\label{eq44}
B^{(k)}_n(f;x)
=\mathbb E
\frac{n!}
     {\big(n-|k|
      \big)!}\,\Delta^kf
\big(n^{-1}\eta_{n-k}(x)
\big);
\end{align}
for the definition of a random vector $\eta_{n}(x)$ see \S\,\ref{Sec3}.

{\bf 5.} 
For the proof by induction of the relations~(\ref{eq4}) and~(\ref{eq44}), let us note that the basis of induction ($k=0$) coincides with the definition of the  polynomial $B_n(f;x)$; as for the inductive step, for the sake of clarity of the following bulky computations let us first  differentiate once the function~$B_n(f;x)$. we have, 
\begin{align*}
\partial_{x_1}B_n(f;x)
&:=\sum_{\substack{0\le j_1{+}j_2{+}\cdots{+}j_d\le n,\\j_1,j_2,\dots j_d\ge 0}}
f\left(\frac{j_1}n,
       \frac{j_2}n,\dots,
       \frac{j_d}n
 \right)\\\noalign{\vskip4pt}
&\qquad\qquad\times
C^j_n\partial_{x_1}x_1^{j_1}
                   x_2^{j_2}\cdots
                   x_d^{j_d}(1{-}x_1{-}x_2{-}\cdots{-}x_d)^{n{-}|j|}\\\noalign{\vskip4pt}
&=\sum_{\substack{0\le j_1{+}j_2{+}\cdots{+}j_d\le n,\\j_1>0,\,
                  j_2,\dots,j_d\ge 0}}
f\bigg(\frac{j_1}n,
       \frac{j_2}n,\dots,
       \frac{j_d}n
 \bigg)\\\noalign{\vskip4pt}
&\qquad\qquad\times
C^j_nj_1x_1^{j_1{-}1}x_2^{j_2}\cdots x_d^{j_d}
            (1{-}x_1{-}x_2{-}\cdots{-}x_d)^{n{-}|j|}\\\noalign{\vskip4pt}
&\qquad-\sum_{\substack{0\le j_1{+}j_2{+}\cdots{+}j_d<n,\\j_1,j_2,\dots,j_d\ge 0}}
f\bigg(\frac{j_1}n,
       \frac{j_2}n,\dots,
       \frac{j_d}n
 \bigg)\\\noalign{\vskip4pt}
&\qquad\qquad\times
C^j_n\big(n{-}|j|
     \big)x_1^{j_1}x_2^{j_2}\cdots x_d^{j_d}
     (1{-}x_1{-}x_2{-}\cdots{-}x_d)^{n{-}|j|{-}1}\\
&=\sum_{\substack{0\le j_1{+}j_2{+}\cdots{+}j_d\le n{-}1,\\j_1,j_2,\dots,j_d\ge 0}}
f\bigg(\frac{j_1{+}1}n,
       \frac{j_2}n,\dots,
       \frac{j_d}n
 \bigg)\\\noalign{\vskip4pt}
&\qquad\qquad\times
C^{j_1{+}1,\dots,j_d}_n(j_1{+}1)
          x_1^{j_1}
          x_2^{j_2}\cdots
          x_d^{j_d}(1{-}x_1{-}x_2{-}\cdots{-}x_d)^{n{-}|j|}\\\noalign{\vskip4pt}
&\qquad-\sum_{\substack{0\le j_1{+}j_2{+}\cdots{+}j_d<n,\\j_1,j_2,\dots,j_d\ge 0}}
f\bigg(\frac{j_1}n,
       \frac{j_2}n,\dots,
       \frac{j_d}n
 \bigg)\\\noalign{\vskip4pt}
&\qquad\qquad\times
C^{j_1,\dots, j_d}_n\big(n{-}|j|
                    \big)x_1^{j_1}
                         x_2^{j_2}\cdots
                         x_d^{j_d}
(1{-}x_1{-}x_2{-}\cdots{-}x_d)^{n{-}1{-}|j|}\\\noalign{\vskip4pt}
&=\sum_{\substack{0\le j_1{+}j_2{+}\cdots{+}j_d<n,\\j_1,j_2,\dots,j_d\ge 0}}
n\Delta_{(x_1)}f\bigg(\frac{j_1}n,
                      \frac{j_2}n,\dots,
                      \frac{j_d}n
                \bigg)\\\noalign{\vskip4pt}
&\qquad\qquad\times
C^{j_1,\dots,j_d}_{n{-}1}x_1^{j_1}
                       x_2^{j_2}\cdots
                       x_d^{j_d}(1{-}x_1{-}x_2{-}\cdots{-}x_d)^{n{-}1{-}|j|},
\end{align*}
due to the identities 
\begin{gather*}
C^{j_1+1,j_2,\dots,j_d}_n(j_1+1)=
\frac{n!(j_1+1)}
     {(j_1+1)!j_2!\dots}=n
\frac{(n-1)!}
     {j_1!j_2!\cdots j_d!}=n
C^{j_1,j_2,\dots,j_d}_{n-1},\\
\noalign{\vskip7pt}
C^{j_1,j_2,\dots,j_d}_n(n-j_1-\cdots-j_d)=
\frac{n!(n-j_1-\cdots-j_d)}
     {j_1!j_2!\cdots\big(n-|j|
                   \big)!}=n
C^{j_1,j_2,\dots,j_d}_{n-1}\,.
\end{gather*}

{\bf 6.} 
The ``full'' induction step: under the  assumption that the formula (\ref{eq4}) holds true for some multi-index $k=(k_1, \dots, k_d)$, let us differentiate this expression again, say, with respect to $x_1$. We get, 
\begin{align*}
\partial_{x_1}B^{(k)}_n(f;x)
&=\partial_{x_1}
  \sum_{0\le j_1+j_2+\cdots+j_d\le n-|k|}
  \frac{n!}
       {\big(n-|k|
        \big)!}\,\Delta^kf
\bigg(\frac{j_1}n,
      \frac{j_2}n,\dots,
      \frac{j_d}n
\bigg)\\\noalign{\vskip4pt}
&\quad\times
C_{n{-}|k|}^{j_1,j_2,\dots,j_d}x_1^{j_1}
                               x_2^{j_2}\cdots\big(1{-}\|x\|
                                             \big)^{n{-}|k|{-}|j|}\\\noalign{\vskip4pt}
&=\sum_{\substack{0\le j_1+j_2+\cdots+j_d\le n,\\j_1>0,j_2,\dots,j_d\ge 0}}
  \frac{n!}
       {\big(n-|k|
        \big)!}\,\Delta^kf
\bigg(\frac{j_1}n,
      \frac{j_2}n,\dots,
      \frac{j_d}n
\bigg)\\\noalign{\vskip4pt}
&\quad\times
C^j_{n{-}|k|}j_1x_1^{j_1{-}1}
                x_2^{j_2}\cdots
                x_d^{j_d}(1{-}x_1{-}
                              x_2{-}\cdots{-}
                              x_d)^{n{-}|k|{-}|j|}\\\noalign{\vskip4pt}
&\ \ -\sum_{\substack{0\le j_1+j_2+\cdots+j_d<n,\\j_1,j_2,\dots,j_d\ge 0}}
      \frac{n!}
           {\big(n-|k|
            \big)!}\,\Delta^kf
      \bigg(\frac{j_1}n,
            \frac{j_2}n,\dots,
            \frac{j_d}n
      \bigg)\\\noalign{\vskip4pt}
&\quad\times
C^j_{n{-}|k|}\big(n{-}|k|{-}|j|
             \big)x_1^{j_1}
                  x_2^{j_2}{\cdots}
                  x_d^{j_d}
(1{-}x_1{-}x_2{-}{\cdots}{-}x_d)^{n{-}|k|{-}|j|{-}1}\\
&=\sum_{\substack{0\le j_1+j_2+\cdots+j_d\le n-1,\\j_1,j_2,\dots j_d\ge 0}}
  \frac{n!}
       {\big(n-|k|
        \big)!}\,\Delta^kf
\bigg(\frac{j_1+1}{n},
      \frac{j_2}{n},\dots,
      \frac{j_d}{n}
\bigg)\\
&\quad\times
C^{j_1{+}1,\dots,j_d}_{n{-}|k|}(j_1{+}1)x_1^{j_1}
                                        x_2^{j_2}{\cdots}
                                        x_d^{j_d}
(1{-}x_1{-}x_2{-}{\cdots}{-}x_d)^{n{-}|k|{-}|j|}\\
&\ \ -\sum_{\substack{0\le j_1+j_2+\cdots+j_d<n,\\j_1,j_2,\dots j_d\ge 0}}
\frac{n!}
     {\big(n-|k|
      \big)!}\,\Delta^kf
\bigg(\frac{j_1}n,
      \frac{j_2}n,\dots,
      \frac{j_d}n
\bigg)\\
&\quad\times
C^{j_1,\dots,j_d}_{n{-}|k|}
\big(n{-}|k|{-}|j|
\big)x_1^{j_1}x_2^{j_2}{\cdots}
              x_d^{j_d}
(1{-}x_1{-}x_2{-}{\cdots}{-}x_d)^{n{-}|k|{-}1{-}|j|}\\
&=\sum_{\substack{0\le j_1+j_2+\cdots+j_d<n,\\j_1,j_2,\dots j_d\ge 0}}
\big(n-|k|
\big)
\frac{n!}
     {\big(n-|k|
      \big)!}\,\Delta_{(x_1)}
               \Delta^kf
\bigg(\frac{j_1}n,
      \frac{j_2}n,\dots,
      \frac{j_d}n
\bigg)\\
&\quad\times
C^{j_1,\dots,j_d}_{n{-}|k|{-}1}x_1^{j_1}
                               x_2^{j_2}\cdots
                               x_d^{j_d}
(1{-}x_1{-}x_2{-}\cdots{-}x_d)^{n{-}|k|{-}1{-}|j|}\\
&=\sum_{\substack{0\le j_1+j_2+\cdots+j_d<n,\\j_1,j_2,\dots j_d\ge 0}}
  \frac{n!}
       {\big(n-|k|-1
        \big)!}\,\Delta_{(x_1)}
                 \Delta^kf
\bigg(\frac{j_1}n,
      \frac{j_2}n,\dots,
      \frac{j_d}n
\bigg)\\
&\quad\times
C^{j_1,\dots,j_d}_{n{-}|k|{-}1}x_1^{j_1}x_2^{j_2}\cdots x_d^{j_d}
(1{-}x_1{-}x_2{-}\cdots{-}x_d)^{n{-}|k|{-}1{-}|j|},
\end{align*}
\vskip-2pt\noindent
by virtue of the identities 
\vskip-7pt\noindent
\begin{gather*}
\aligned
C^{j_1+1,j_2,\dots,j_d}_{n-|k|}(j_1+1)
&=\frac{(n-|k|)!(j_1+1)}
       {(j_1+1)!j_2!\cdots\big(n-|k|-|j|-1
                         \big)!}\\
&=\big(n-|k|
  \big)\frac{(n-|k|-1)!}
            {j_1!j_2!\cdots\big(n-|k|-1-|j|
                          \big)!}\\
&=\big(n-|k|
  \big)C^{j_1,j_2,\dots,j_d}_{n-|k|-1},
\endaligned\\
\aligned
C^{j_1,j_2,\dots,j_d}_{n{-}|k|}\big(n{-}|k|{-}j_1{-}\cdots{-}j_d
                             \big)
&=\frac{\big(n-|k|
        \big)!\big(n-|k|-|j|
              \big)}
       {j_1!j_2!\cdots,j_d!\big(n-|k|-|j|
                         \big)!}\\
&=\big(n-|k|
 \big)C^{j_1,j_2,\dots,j_d}_{n-|k|-1}.
\endaligned
\end{gather*}
Hence, by induction the formula~(\ref{eq4}) holds true along with~(\ref{eq44}).

{\bf 7.} Again having in mind how Bernstein's method for the Theorem \ref{thm1-en} was applied in the proof of the Theorem \ref{thm4-en} for the  polynomials $B_n$, and similarly to what was done earlier for $\wtb_n$, let us recall that for $f\in C^{k}(\mathbb R^d)$ the normalized finite differences are uniformly close to the corresponding partial derivatives, namely,
$$
\frac{n!}
     {\big(n-|k|
      \big)!}\,\Delta^kf(x)\rightrightarrows
\frac{\partial^{|k|}f}
     {\partial x^k}(x),
\ \ \mbox{for}\ \ x\in\mathbb S^d,
\ \ n\to\infty
$$
(cf. the Lemma \ref{lem1}). Hence statement of the Theorem for the polynomials $B_n$ follows from the law of large numbers~(\ref{llneta}) and from the representation~(\ref{eq44}).

The Theorem~\ref{thm4-en} is proved
\goodbreak

\subsection*{Acknowledgement}

The authors are grateful to the anonymous referee for useful advice.


\begin{thebibliography}{99}

\bibitem{ia}
{\sl Abel,~U. and Ivan,~M.}
Asymptotic expansion of the~multivariate Bernstein polynomials
on a~simplex~// %
{\it Approx. Theory Appl.~$N.\,S.$}. 2000. V.\,16, N\,3. P.\,85--93.



\bibitem{B2}
{\sl Bernstein, S.~N.}
Complement \`a~l'article de E. Voronovskaja 
``Determination de la forme asymptotique de l'approximation des fonctions par les polyn\^omes de S. Bernstein''~//
{\it Dokl. Akad. Nauk SSSR, Ser.~A (Russian)}, 1932. V.\,4. P.\,86--92.

\bibitem{bern}
{\sl Bernstein, S.~N.}
Communications de la Soci\'et\'e math\'ematique de Kharkow. 2-\'ee s\'erie, 1912, V.\, 13, N\,1, P.\,1--2.
(See also http://topo.math.u-psud.fr/$\sim$slc/TeX/\break lcs/legacy\_math/bernstein/bernstein\_doyle.pdf)


\bibitem{Hlo}
{\sl Chlodovsky, ~I.~N.} 
On some properties of SN Bernstein polynomials~//
Moscow, Proc. 1st All-Union congress of mathematics (Russian), (Kharkov, 1930). P.\,22.



\bibitem{deriv}
{\sl Floater, M.~S.}
On~the convergence of derivatives of Bernstein approximation~// %
{\it J.~Approx. Theory}. 2005. V.\,134, N\,1. P.\,130--135.

\bibitem{Kry}
{\sl Krylov, N.~V.}
{\it Introduction to the Theory of Random Processes}
AMS, Providence, Rhode Island, 2002.


\bibitem{L1937}
{\sl Lorentz, G.~G.}
Zur Theorie der Polynome von S.~Bernstein~// %
{\it Mathematical Collections} 1937. V.\,2(44), No.\,3. P.\,543--556.  

\bibitem{L}
{\sl Lorentz, G.~G.}
{\it Bernstein Polynomials}~/~2nd ed.
New~York: AMS Chelsea Publishing Co., 1986.

\bibitem{review}
{\sl Pop, O.~T.}
Voronovskaja-type theorem for certain GBS operators~// %
{\it Glas. Mat. Ser~III}. 2008. V.\,43(63), N\,1. P.\,179--194.

\bibitem{pop35}
{\sl Popoviciu, T.}
Sur l'approximation des fonctions convexes d'ordre sup\-\'erieur~// %
{\it Mathematica}. 1935. V.\,10. P.\,49--54.


\bibitem{ash}
{\sl Shiryaev, A.~N}
{\it Probability}
New~York, Springer, 2nd edition, 1995.


\bibitem{tel}
{\sl Telyakovskii, S.~A.}
On the rate of approximation of functions by the Bernstein polynomials~// %
{\it Proc. Inst. Math. Mech.}, 2009. V.\,264, suppl.\,1. P.\,177--184.

\bibitem{Ti-Sh}

{\sl Tikhonov, I.~V., Sherstyukov, V.~B.}
The module function approximation by Bernstein polynomials~//
{\it Bulletin of Chelyabinsk State University (Russian)} 2012. No. 26(15). P.\,6-40.


\bibitem{VV2012}
{\sl Veretennikov, A.~Yu., Veretennnikova, E.~V.}
On convergence of partial derivatives of multivariate Bernstein polynomials~//
{\it Mathematics, Informatics and Physics in Science and Education}~/
Collection of scientific papers for 140th anniversary of MSPU (Russian). Moscow, ``Prometei'', 2012. P.\,39--42.

\bibitem{rate}
{\sl Voronovskaya, E.~V.}
Determination of the asymptotic form of approximation of functions by the Bernstein polynomials~//
{\it Dokl. Akad. Nauk SSSR, Ser.~A (Russian)}, 1932. V.\,4. P.\,74--85.

\bibitem{Zor2}
{\sl Zorich, V.~A.}
{\it Mathematical Analysis II}. 
Berlin, Springer, 2004.






\end{thebibliography}
\end{document}